\documentclass[11pt]{amsart}
\usepackage{mathrsfs}
\usepackage{}
\usepackage{xypic}
\usepackage{amsfonts}
\usepackage{amssymb}
\usepackage{bbm}

\newtheorem{theo}{Theorem}[section]

\newtheorem{lem}[theo]{Lemma}
\newtheorem{exa}[theo]{Example}
\newtheorem{prop}[theo]{Proposition}
\newtheorem{defi}[theo]{Definition}
\newtheorem{rem}[theo]{Remark}
\numberwithin{equation}{section}

\newcommand{\be}{\begin{equation}}
\newcommand{\ee}{\end{equation}}
\newcommand{\bes}{\begin{eqnarray}}
\newcommand{\ees}{\end{eqnarray}}
\newcommand{\bess}{\begin{eqnarray*}}
\newcommand{\eess}{\end{eqnarray*}}

\newcommand{\bali}{\begin{align}}
\newcommand{\eali}{\end{align}}
\newcommand{\epf}{\hfill \rule {2mm}{3mm}\vskip10pt}

\newcommand{\pf}{\noindent{\bf Proof.}\quad}

\begin{document}
\title[Green Rings of Pointed Rank One Hopf algebras]{Green Rings of Pointed Rank One Hopf algebras of Nilpotent Type}
\author{Zhihua Wang}
\address{School of Mathematical Science, Yangzhou University,
Yangzhou 225002, China}
\email{mailzhihua@gmail.com}
\author{Libin Li}
\address{School of Mathematical Science, Yangzhou University,
Yangzhou 225002, China}
\email{lbli@yzu.edu.cn}
\author{Yinhuo Zhang}
\address{Department of Mathematics and Statistics, University of Hasselt, Universitaire Campus, 3590 Diepeenbeek,Belgium }
\email{yinhuo.zhang@uhasselt.be}
\date{}
\subjclass[2000]{16W30, 19A22} \keywords{Green ring, indecomposable
module, symmetric ring, Jacobson radical, group-like algebra, bi-Frobenius algebra}

\begin{abstract} Let $H$ be a finite dimensional
pointed rank one Hopf algebra of nilpotent type. We first determine all finite
dimensional indecomposable $H$-modules up to isomorphism, and then establish the Clebsch-Gordan formulas for the decompositions of the tensor
products of indecomposable $H$-modules by virtue of almost split sequences. The Green ring $r(H)$ of $H$ will be presented in terms of generators and relations. It turns out that the Green ring $r(H)$ is commutative and is generated by one variable over the Grothendieck ring $G_0(H)$ of $H$ modulo one relation.  Moreover, $r(H)$ is Frobenius and symmetric with dual bases associated to almost split sequences, and its Jacobson radical is a principal ideal. Finally, we show that the stable Green ring, the Green ring of the stable module category, is isomorphic to the quotient ring of $r(H)$ modulo the principal ideal generated by the projective cover of the trivial module. It turns out that the complexified stable Green algebra is a group-like algebra and hence a bi-Frobenius algebra.
\end{abstract}

\maketitle

\section{\bf Introduction }

The study of Green rings (or representation rings) of non-semisimple  Hopf algebras has been recently revived. In \cite{CVZ} and \cite{LZ}, the authors studied respectively the Green rings of the Taft algebras and the generalized Taft algebras  based on the Cibils' work \cite{Ci} on the decomposition of tensor products of indecomposable modules. The nilpotent elements of the Green rings of the (generalized) Taft algebras were completely determined in terms of linear combinations of projective indecomposable modules.  In \cite{HOYZ}, the Green rings of  pointed tensor categories of finite type were investigated where the quiver techniques were applied to study the comodules over a graded coquasi-Hopf algebra. An interesting phenomenon appears when one compares the Green ring of a generalized Taft algebra $T$ and the Green ring of its dual $T^*$. Although the module category of $T$ and the comodule category of $T$ are not monoidal equivalent, but the Green rings of the two are isomorphic. This means that the Green ring invariant of a Hopf algebra can not solely determine the module (or comodule) category of the Hopf algebra. For Hopf algebras of infinite type, the Green rings are usually not finitely generated. In this case, the Green rings are difficult to study. For example, the Green ring of the Drinfeld double of the Sweedler Hopf algebra \cite{Ch} and the Green rings of the rank 2 Taft algebras \cite{LH} are not finitely generated although they can be computable.

In this paper, we study the Green rings of a large class of finite dimensional pointed Hopf algebras of finite representation type. The pointed Hopf algebras in this class are of rank one, which are  either of nilpotent type or of non-nilpotent type.
Let $A$ be a Hopf algebra over a field $\mathbbm{k}$ and $A_0\subseteq
A_1\subseteq A_2\subseteq\cdots$ the coradical filtration of $A$.
Assume that the coradical $A_0$ is a Hopf subalgebra of $A$. Then
each $A_i$ is a free $A_0$-module. Consider $\mathbbm{k}$ as the trivial right $A_0$-module. If $A$ is generated as an algebra by $A_1$ and $\dim_{\mathbbm{k}}(\mathbbm{k}\otimes_{A_0}A_1)=n+1$,
then $A$ is called a Hopf algebra of rank $n$ \cite{KR, WYC}. The (generalized) Taft algebras and the half quantum group \cite{Gu} are all finite dimensional pointed Hopf algebras of rank one.
Any finite dimensional pointed Hopf algebra of rank one can be constructed from a  group datum. So the classification of finite dimensional pointed Hopf algebras of rank one over an algebraically closed field (of characteristic $0$) can be derived from the classification of the group data, see \cite{KR, CHYZ}.

The Green ring $r(H)$ of a finite dimensional pointed Hopf algebra $H$ of rank one is closely related to the Grothedieck ring $G_0(\mathbbm{k}G))$ of the group algebra of $G$, where $G$ is the group of group-like elements of $H$. In the case of nilpotent type, $r(H)$ is isomorphic to some quotient ring of the polynomial ring $G_0(\mathbbm{k}G)[z]$ modulo one interesting relation presented by a Dickson polynomial multiplied with the almost split sequence ending with the trivial module. The Green ring $r(H)$ possesses a natural Frobenius structure and is a symmetric ring. The Jacobson radical of $r(H)$ is a principal ideal generated by a special element presented by a linear combination of projective modules. Moreover, we will show that the Green ring $r_{st}(H)$ of the stable category of $H$-modules is isomorphic to the quotient ring of $r(H)$ modulo the ideal generated by all projective modules. The complexified Green ring $R_{st}(H)$ is both a group-like algebra and a bi-Frobenius algebra. However, if $H$ is of non-nilpotent type, the Green ring $r(H)$ is not nicely generated over the Grothendieck ring, but it possesses the Frobenius structure as well. We will present the structure of the Green ring in a separate paper \cite{WLZ2}. This paper is organized as follows.

In Section 2, we first recall the construction of a finite dimensional pointed Hopf algebra $H$ of rank one from a group datum. If $H$ is of nilpotent type, we show that $H$ is a Nakayama algebra and determine all finite dimensional indecomposable $H$-modules up to isomorphism. In order to obtain the decompositions of tensor products of indecomposable $H$-modules, we study almost split sequences of indecomposable modules in Section 3.  From the uniqueness of an almost split sequence, we can deduce the Clebsch-Gordan formulas for the tensor product of two indecomposable modules. This approach, in fact, works for more general pointed Hopf algebras stemming from quivers.

In Section 4, we use the obtained decomposition formulas to deliver the structure of the Green ring $r(H)$ in terms of generators and relations. More explicitly, the Green ring $r(H)$ is isomorphic to a polynomial ring with one variable over the Grothendieck ring $G_0(H)$ modulo a relation given by a Dickson polynomial multiplied with the almost split sequence ending with the trivial module. As an example, we describe the Green ring of a pointed rank one Hopf algebra $H$ with $G(H)$  being a dihedral group. In this case, the Green ring of $H$ can be explicitly presented with four generators and five relations over $\mathbb{Z}$.

It is well-known that two non-gauge equivalent group algebras may possess the same character ring. This happens to Hopf algebras as well. In Section 4, we shall characterize the group datum, and give a sufficient condition for two non-gauge equivalent rank one Hopf algebras of nilpotent type to share the same Green ring. Taft algebras give such examples. In the end of this section, we shall define an associative symmetric and nonsingular $\mathbb{Z}$-bilinear form on the Green ring $r(H)$. Consequently, $r(H)$ is symmetric with a pair of dual bases associated with the almost split sequences. Thus various ring-theoretical properties can be derived from the fact that $r(H)$ is both Frobenius and symmetric.

In Section 5, we study the Jacobson radical $J(r(H))$ of the Green ring $r(H)$. To this end, we consider the complexified Green ring $R(H):=\mathbb{C}\otimes_{\mathbb{Z}}r(H)$ over the field $\mathbb{C}$. We can determine the dimension of the Jacobson radical of $R(H)$ by calculating the number of simple modules over $R(H)$.  The rank of the Jacobson radical of $r(H)$ is in fact equal to the dimension of the Jacobson radical of $R(H)$. It turns out that the Jacobson radical of $r(H)$ is a principal ideal generated by a special element represented by a linear combination of projective modules (Theorem \ref{th4.5}). This explains the reason why those nilpotent elements of the Green rings of the generalized Taft algebras are of a special linear combination form of projective modules \cite{LZ}.

In Section 6, we study the Green ring of the stable category of left $H$-modules. We show that the Green ring of the stable category is the quotient ring of $r(H)$ modulo the principal ideal generated by the projective cover of the trivial module. Furthermore, the complexified stable Green algebra possesses a group-like algebra structure and a bi-Frobenius algebra structure.  As a consequence, many interesting properties of group-like algebras and bi-Frobenius algebras (see \cite{Doi1, Doi2, Doi3, DT}) hold for the complexified stable Green algebra.

Throughout, we work over an algebraically closed field $\mathbbm{k}$
of characteristic zero. Unless otherwise stated, all algebras, Hopf
algebras and modules are defined over $\mathbbm{k}$; all modules are
left modules and finite dimensional; all maps are
$\mathbbm{k}$-linear; $\dim$, $\otimes$ and $\textrm{Hom}$ stand for
$\dim_{\mathbbm{k}}$, $\otimes_{\mathbbm{k}}$ and
$\textrm{Hom}_{\mathbbm{k}}$, respectively.  In the sequel, we will denote by $H$  a finite dimension pointed rank one Hopf algebra of nilpotent type, and let $A$ stand for a non-specified Hopf algebra. For the theory of Hopf algebras, we refer to \cite{Mon, Swe}.

\section{\bf Representations of pointed rank one Hopf algebras of nilpotent type}
\setcounter{equation}{0}
\def\kk{\mathbbm{k}}
\def\Hd{H_{\mathcal {D}}}

In this section, we first recall the construction of a finite
dimensional pointed Hopf algebra of rank one from a group datum $\mathcal
{D}$.  A quadruple $\mathcal {D}=(G,\chi,g,\mu)$ is
called a \emph{group datum} if $G$ is a finite group, $g$ an element in
the center of $G$, $\chi$  a $\mathbbm{k}$-linear character of $G$,
and $\mu\in\mathbbm{k}$ subject to $\chi^n=1$ or $\mu(g^n-1)=0$,
where $n$ is the order of $\chi(g)$. If $\mu(g^n-1)=0,$ then the
group datum $\mathcal {D}$ is said to be of \emph{nilpotent type}.
Otherwise, it is of \emph{non-nilpotent type}, see \cite{KR}.

Given a group datum $\mathcal {D}=(G,\chi,g,\mu)$, we let
$H_{\mathcal {D}}$ be an associative algebra generated by $y$ and all $h$ in $G$ such that $\kk G$ is a subalgebra of $H_{\mathcal {D}}$ and
$$y^n=\mu(g^n-1),\ yh=\chi(h)hy,\ \textrm{for}\ h\in G.$$

The algebra $\Hd$ is finite dimensional with a canonical PBW
$\mathbbm{k}$-basis $\{y^ih|h\in G,\ 0\leq i\leq n-1\}$.
Thus $\dim H_{\mathcal {D}}=n|G|$, where $|G|$ is the order of $G$.

In fact, $H_{\mathcal {D}}$ is endowed with a Hopf algebra
structure. The comultiplication $\bigtriangleup$, the counit
$\varepsilon$, and the antipode $S$ are given respectively by
$$\bigtriangleup(y)=y\otimes g+1\otimes y,\ \varepsilon(y)=0,\ S(y)=-yg^{-1},$$
$$\bigtriangleup(h)=h\otimes h,\ \varepsilon(h)=1,\ S(h)=h^{-1},$$
for all $h\in G.$

It is easy to see that $G$ is the group of group-like elements of
$H_{\mathcal {D}}$ and $\Hd$ is a pointed Hopf algebra of rank one \cite{KR}. If the group datum $\mathcal {D}$ is of nilpotent type, $\Hd$ is said to be a \emph{pointed rank one Hopf algebra of nilpotent type}. Otherwise, it is of \emph{non-nilpotent type}.

The linear
character $\chi$ induces an automorphism $\sigma$ of $\kk G$:
$$\sigma(a)=\Sigma\chi(a_1)a_2,$$
where $a\in \kk G$ and $\bigtriangleup(a)=\sum a_1\otimes a_2.$
In this case, $y^ja=\sigma^j(a)y^j$ holds in $H_{\mathcal {D}}$, for $j\geq 0$. The following proposition gives a classification of Hopf algebras $H_{\mathcal {D}}$ with group data $\mathcal {D}$, see \cite[Theorem 1 (c)]{KR}, or \cite[Theorem 5.9]{CHYZ}.

\begin{prop}\label{prop2}
Let $\mathcal {D}=(G,\chi,g,\mu)$ and $\mathcal
{D'}=(G',\chi',g',\mu')$ be two group data. Then the Hopf algebras $H_\mathcal {D}$ and $H_\mathcal {D'}$ are isomorphic if and only if
there is a group isomorphism $f:G\rightarrow G'$ such that
$f(g)=g'$, $\chi=\chi'\circ f$ and $\beta\mu'(g'^n-1)=\mu(g'^n-1)$
for some non-zero $\beta\in\mathbbm{k}$, where $n$ is the order of
$\chi(g)$.
\end{prop}

In case the group datum $\mathcal {D}=(G,\chi,g,\mu)$ is of
nilpotent type, namely, $\mu(g^n-1)=0$, where $n$ is the order of $\chi(g)$, then it is either $\mu=0$ or $g^n-1=0$. In both cases, Proposition \ref{prop2} implies that the Hopf algebras constructed from $(G,\chi,g,\mu)$ and $(G,\chi,g,0)$ respectively are isomorphic. Because of this fact, we may assume that $\mu=0$ for any group datum $\mathcal {D}=(G,\chi,g,\mu)$ of nilpotent type.

From now on, we only deal with finite dimensional pointed rank one Hopf algebras of nilpotent type. We fix the following notations throughout the paper.  Let $H:=H_\mathcal {D}$, where  $\mathcal{D}=(G,\chi,g,0)$. Assume that $n\geq2$ is the order of $q:=\chi(g)$.
Since $q$ is a primitive $n$-th root of unity and $n\geq2$, we have that
$g\neq 1$ and $\chi\neq\varepsilon$. Note that $y^n=0$ and the
quotient algebra $H/(y)\cong\mathbbm{k}G$. Hence, the
Jacobson radical $J$ of $H$ is generated  by the
nilpotent element $y$. Let $e=\frac{1}{|G|}\sum_{h\in G}h$. Then
$ey^{n-1}\in\int_H^l$ and $y^{n-1}e\in\int_H^r$. It follows from
\cite[Proposition 2.5]{Lo} that $H$ is neither unimodular nor
symmetric.

Since $J=(y)$ and $H/J\cong\mathbbm{k}G$, an $H$-module $V$ is simple if and only if $yV=0$ and $V$ restricts to a simple  $\mathbbm{k}G$-module. Thus a complete set of non-isomorphic simple $\mathbbm{k}G$-modules forms a complete set of non-isomorphic simple $H$-modules. In the sequel,  we fix such a complete set $\{V_1,V_2,\cdots,V_m\}$ of non-isomorphic simple $\mathbbm{k}G$ (and $H$)-modules.


\begin{rem}\label{rem2.1}
The fact that each simple $H$-module $V_i$ restricts to a simple $\mathbbm{k}G$-module implies that there exists a primitive idempotent $e_i$ of $\mathbbm{k}G$ such that $V_i\cong\mathbbm{k}Ge_i$. Since $e_i$ is also a primitive idempotent of $H$ and $He_i/JHe_i\cong V_i,$ we see that $He_i$
is the projective cover of the simple $H$-module $V_i$, for $1\leq i\leq
m$.
\end{rem}

Now let $V_{\chi}$ and $V_{\chi^{-1}}$ be the two (one dimensional) simple
$\mathbbm{k}G$-modules corresponding to the linear characters $\chi$
and $\chi^{-1}$ respectively. For any simple $\mathbbm{k}G$-module
$V_i$, $1\leq i\leq m$, the tensor product $V_{\chi^{-1}}\otimes V_i\cong V_i\otimes
V_{\chi^{-1}}$ is simple as well. Hence there is a unique permutation $\tau$ of the index set $\{1,2,\cdots,m\}$ such that
 $$V_{\chi^{-1}}\otimes V_i\cong V_i\otimes
V_{\chi^{-1}}\cong V_{\tau(i)}.$$ The inverse of $\tau$ is determined by $$V_{\chi}\otimes V_i\cong V_i\otimes V_{\chi}\cong V_{\tau^{-1}(i)}.$$

\begin{lem}\label{lem3.1}
For any $1\leq i\leq m$ and $t\in\mathbb{Z}$, there is a bijective map $\widetilde{\sigma}_{i,t}$ from $V_i$ to $V_{\tau^t(i)}$ such that
$\widetilde{\sigma}_{i,t}(av)=\sigma^t(a)\widetilde{\sigma}_{i,t}(v)$,
for any $a\in \mathbbm{k}G$ and $v\in V_i$.
\end{lem}
\pf For a fixed non-zero element $u\in V_{\chi^{-t}}$, the map
$V_i\rightarrow V_{i}\otimes V_{\chi^{-t}},\ v\mapsto v\otimes u$
composed with the isomorphism $V_i\otimes V_{\chi^{-t}}\cong
V_{\tau^t(i)}$ gives the desired bijective map.
\epf

Let $x$ be a variable and $V$ a $\kk G$-module. For any $k\in \mathbb{N}$,
consider $x^kV$ as a vector space defined by
$x^ku+x^kv=x^k(u+v)$ and $\lambda(x^ku)=x^k(\lambda u)$, for $u,v\in V, \lambda\in\mathbbm{k}.$
Then $x^kV$ becomes a $\kk G$-module defined by
\be\label{equ3.3}h(x^kv)=\chi^{-k}(h)x^khv,\ee for any $h\in G$ and $v\in V$.  We have the following lemma.

\begin{lem}\label{lem3.2}
For any simple $\kk G$-module $V_i$ and $k\in \mathbb{N}$, we have $\kk G$-module isomorphisms $x^kV_i\cong V_{\chi^{-k}}\otimes V_i\cong V_{\tau^k(i)}$.
\end{lem}

For any $1\leq i\leq m$ and $1\leq j\leq n$, the direct sum $$M(i,j):=V_i\oplus xV_i\oplus\cdots\oplus x^{j-1}V_i$$
 is a $\kk G$-module, where each summand is a simple $\kk G$-module defined by (\ref{equ3.3}). We define an action of $y$ on $M(i,j)$ as follows:
\be\label{equ3.4}y(x^kv)=\begin{cases}
x^{k+1}v, & 0\leq k\leq j-2,\\
0, & k=j-1,
\end{cases}\ee
for any $v\in V_i$. Then $M(i,j)$ becomes an $H$-module with $\dim M(i,j)=j\dim(V_i)$. Moreover, it is easy to see that $M(i,1)\cong V_i$ and $M(i,n)\cong He_i$. If $\mathbf{v}_i$ is a basis of $V_i$, then $\{x^kv\mid0\leq k\leq j-1,v\in\mathbf{v}_i\}$ forms a basis of $M(i,j)$. In particular, $\{1,x,\cdots,x^{j-1}\}$ forms a basis of $M(1,j)$, where we identify  $x^k1$ with $x^k$, $1\in \kk$.

\begin{theo}\label{th2.2} Let $JM(i,j)$ and $P(M(i,j))$ be the radical and the projective cover of $M(i,j)$ respectively, for $1\leq i\leq m$, $1\leq j\leq n$. Then
\begin{enumerate}
\item[(1)] $JM(i,1)=0$ and $JM(i,j)\cong M(\tau(i),j-1)$, for $2\leq j\leq n$.
\item[(2)] $\textrm{soc}M(i,j)\cong V_{\tau^{j-1}(i)}$, $M(i,j)/JM(i,j)\cong V_i$ and
$P(M(i,j))\cong M(i,n)$.
\item[(3)] $M(i,j)$ is indecomposable and uniserial. $H$ is a Nakayama algebra, and hence it is of finite representation type.
\item[(4)] $M(i,j)\cong M(k,l)$ if and only if $i=k$ and $j=l$. Moreover, the set
$\{M(i,j)\mid1\leq i\leq m,\ 1\leq j\leq n\}$ forms a complete set of
finite dimensional indecomposable $H$-modules up to isomorphism.
\end{enumerate}
\end{theo}
\pf (1). Since $M(i,1)\cong V_i$ is simple, we have $JM(i,1)=0$. For $2\leq j\leq n$, as $J=(y)$, we have $JM(i,j)=\oplus_{k=1}^{j-1}x^kV_i$ as a vector space.
Define a $\kk$-linear map $$\phi:JM(i,j)\rightarrow M(\tau(i),j-1),\ x^kv\mapsto x^{k-1}\widetilde{\sigma}_{i,1}(v),$$ for any $v\in V_i$ and $1\leq k\leq j-1.$ By Lemma \ref{lem3.1}, it is straightforward to check that the map $\phi$ is an $H$-module isomorphism.

(2) It is easy to  see that $\textrm{soc}M(i,j)=\{u\in M(i,j)\mid yu=0\}=x^{j-1}V_i\cong V_{\tau^{j-1}(i)}$ as $H$-modules (because of  $J=(y)$).
Since $JM(i,j)=\oplus_{k=1}^{j-1}x^kV_i$, we have
$M(i,j)/JM(i,j)\cong V_i$. The following isomorphism
$$P(M(i,j))\cong P(M(i,j)/JM(i,j))\cong P(V_i)\cong
M(i,n)$$ follows from Remark \ref{rem2.1}.

(3) That $M(i,j)$ is indecomposable follows from the fact that $\textrm{soc}M(i,j)$ is simple. The proof of $M(i,j)$ to be uniserial is similar to the one of  \cite[Proposition 3]{KR}. Denote by $N_{i,l}$ the
submodule of $M(i,j)$:
$$N_{i,l}=x^lV_i\oplus\cdots\oplus x^{j-1}V_i,\ \textrm{for}\ 0\leq l\leq j-1.$$
Suppose that $N$ is a non-zero submodule of
$M(i,j)$. Then there exists a largest $l$ such that $N\subseteq
N_{i,l}$. Since $N_{i,l+1}$ is a maximal submodule of $N_{i,l}$, we
conclude that $N+N_{i,l+1}=N_{i,l}$. However, $N_{i,l+1}=JN_{i,l}$,
whence $N=N_{i,l}$ by Nakayama's lemma. Thus $M(i,j)$ is uniserial.
Since $H$ is Frobenius (hence self-injective), the indecomposable projective modules $M(i,n)$, $1\leq i\leq m$, are also injective modules.  Hence they are uniserial. It follows that $H$ is a Nakayama algebra.

$(4)$ If $M(i,j)\cong M(k,l),$ then by part $(2)$,
$$V_i\cong M(i,j)/JM(i,j)\cong M(k,l)/JM(k,l)\cong V_k,$$
for $1\leq i,k\leq m$. This implies that $i=k$. Comparing the
dimensions of the vector spaces, we obtain that $j=l$. Since $H$ is a Nakayama algebra, every indecomposable $H$-module $M$ is a
quotient of an indecomposable projective module $M(i,n)$ for
some $1\leq i\leq m$. Thus $M$ is of the form $M(i,j)$, for some $1\leq j\leq n$. \epf

\section{\bf Almost split sequences and Clebsch-Gordan formulas}
\setcounter{equation}{0} Almost split sequences over Nakayama algebras have been much studied,  see, e.g. \cite{ARS}. In this section, we first point out that the Auslander-Reiten translate of the $H$-module $M(i,j)$ is nothing but $M(\tau(i),j)$. We then show that the almost split sequence ending with the non-projective module $M(i,j)$ can be obtained by tensoring  $M(i,j)$ over $\kk$ on the right (or on the left) with the almost split sequence of the trivial $H$-module. This approach using almost split sequences works also for the $\mathscr{H}$-modules, where $\mathscr{H}$ is a Hopf algebra associated to the quiver $A_\infty^\infty$, see \cite{Ci}. But it does not work for more general Hopf algebras, see e.g. \cite{GMS}. After that we use the uniqueness of an almost split sequence to determine the decomposition of the tensor products of indecomposable $H$-modules. For the necessary background of almost split sequences, we refer to the books \cite{ARS, ASS}. To begin with, we need the following proposition and lemma.

\begin{prop}\label{prop3.1}
For any $1\leq i\leq m$ and $1\leq j\leq n$, $V_{i}\otimes M(1,j)\cong M(1,j)\otimes V_{i}\cong M(i,j)$. In particular,
$V_{\chi^{-1}}\otimes M(1,j)\cong M(1,j)\otimes V_{\chi^{-1}}\cong M(\tau(1),j)$.
\end{prop}
\pf The $H$-module isomorphism from $V_i\otimes
M(1,j)$ to $ M(i,j)$ is given by $v\otimes x^k\mapsto x^kv,$ for $0\leq k\leq j-1,\ v\in V_i$. Similarly, let $\lambda_i$ be the scalar such that $gv=\lambda_iv$, for $v\in V_i$. Then the $H$-module isomorphism from $M(1,j)\otimes V_{i}$ to $M(i,j)$ is given by $x^k\otimes v\mapsto \lambda_i^{-k}x^kv$, for $0\leq k\leq j-1,\ v\in V_i$. \epf

\begin{lem}\label{lem3.3}
For any $1\leq i\leq m$ and $1<j<n$, we have the following:
\begin{enumerate}
\item[(1)]There is an injective morphism from $JM(i,j)$ to $M(1,2)\otimes
M(i,j)$.
\item[(2)]There is an injective morphism from $JM(i,j)$ to $
M(i,j)\otimes M(1,2)$.
\end{enumerate}
\end{lem}
\pf We only prove part (1) and the proof of part (2) is similar. Let $\lambda_i$ be a scalar such that $gv=\lambda_iv$, for any $v\in V_i$. Denote by $N$ the
subspace of $M(1,2)\otimes M(i,j)$ spanned by the elements
$$1\otimes x^kv-\zeta_kx\otimes x^{k-1}v,\ \textrm{for\ all}\ v\in V_i,\ 1\leq k\leq j-1,$$
where $\zeta_k=\frac{\lambda_i}{1-q}(\frac{1}{q^{j-1}}-\frac{1}{q^{k-1}})$ and $q=\chi(g)$. Then the action of $H$ on $N$ is stable:
$$h(1\otimes x^kv-\zeta_kx\otimes x^{k-1}v)=\chi^{-k}(h)(1\otimes x^khv-\zeta_kx\otimes
x^{k-1}hv)\in N,$$ for $h\in G$, $1\leq k\leq j-1$,
$$y(1\otimes x^kv-\zeta_kx\otimes x^{k-1}v)=1\otimes x^{k+1}v-\zeta_{k+1}x\otimes x^{k}v\in
N,$$ for $1\leq k\leq j-2$ and $y(1\otimes x^{j-1}v-\zeta_{j-1}x\otimes x^{j-2}v)=0.$
Now we define a $\mathbbm{k}$-linear map $\iota$ from $JM(i,j)=\oplus_{k=1}^{j-1}x^kV_i$ to $N$ as follows:
\be\label{equ3.2}\iota(x^kv)=1\otimes
x^kv-\zeta_kx\otimes
x^{k-1}v\ \textrm{for}\ 1\leq k\leq j-1.\ee It is obvious
that the map $\iota$ is an $H$-module isomorphism. Hence $\iota$ is injective from $JM(i,j)$ to $M(1,2)\otimes M(i,j)$.
\epf

The permutation $\tau$ of the index set $\{1,2,\cdots,m\}$ is related to the Auslander-Reiten translate of $H$-modules as shown in the following.
\begin{prop}\label{lem7.1}
For each non-projective indecomposable module $M(i,j)$, the Auslander-Reiten translate $\textrm{DTr}(M(i,j))$ of $M(i,j)$ is isomorphic to $M(\tau(i),j)$.
\end{prop}
\pf By Theorem \ref{th2.2}, $M(i,j)$ is uniserial and consequently the length of $M(i,j)$ coincides with the radical length of $M(i,j)$ \cite[Proposition 2.1, ChIV]{ARS}, which is exactly the number $j$. It also follows from Theorem \ref{th2.2} that the projective cover of $M(i,j)$ is $M(i,n)$. Thus, by \cite[Proposition 2.6, ChIV]{ARS}, $\textrm{DTr}(M(i,j))\cong JM(i,n)/J^{j+1}M(i,n)\cong M(\tau(i),j)$, as desired.\epf

Recall that $M(1,2)=V_1\oplus xV_1$, and the summand $xV_1$ is isomorphic to $V_{\chi^{-1}}$ as $\kk G$-modules, where the isomorphism is given by $\rho:V_{\chi^{-1}}\rightarrow xV_1$, $\rho(w)=x$, for a fixed non-zero element $w\in V_{\chi^{-1}}$ and $x\in xV_1$.

\begin{lem}\label{lem4.3}
The sequence \be\label{equ4.3}0\rightarrow
V_{\chi^{-1}}\xrightarrow{\alpha}M(1,2)\xrightarrow{\beta}V_1\rightarrow0\ee is an almost split sequence of $H$-modules, where $\alpha=(0,\rho)^\intercal$ and $\beta=(id,0)$.
\end{lem}
\pf It is obvious that the maps $\alpha$ and $\beta$ are both $H$-module morphisms. The short sequence  (\ref{equ4.3}) is exact but not split since $M(1,2)\cong V_1\oplus V_{\chi^{-1}}$ holds just as $\kk G$-modules but not as $H$-modules. Note that $V_1$ is a non-projective brick ($\textrm{End}(V_1)=\kk$). By Proposition \ref{lem7.1}, $\textrm{DTr}(V_1)=\textrm{DTr}(M(1,1))\cong M(\tau(1),1)\cong V_{\chi^{-1}}$. Hence, the sequence  (\ref{equ4.3}) is almost split according to \cite[Corollary 2.4, ChV]{ARS}.\epf

For ech indecomposable $H$-module $M(i,j)$, we consider the following two short exact sequences  obtained  by tensoring $M(i,j)$ over $\kk$  with the almost split sequence (\ref{equ4.3}) on the right and on the left respectively.
\be\label{equ6.5}0\rightarrow M(\tau(i),j)\rightarrow M(1,2)\otimes M(i,j)\rightarrow M(i,j)\rightarrow0,\ee
\be\label{equ6.6}0\rightarrow M(\tau(i),j)\rightarrow M(i,j)\otimes M(1,2)\rightarrow M(i,j)\rightarrow0.\ee

\begin{prop}\label{th6.2} If $M(i,j)$ is non-projective, that is, $j\neq n$, then the short exact sequences of (\ref{equ6.5}) and (\ref{equ6.6}) ending with $M(i,j)$ are both almost split.
\end{prop}
\pf It is enough to show the sequence (\ref{equ6.5}) since the proof of the sequence (\ref{equ6.6}) is similar. Note that if $j=1$, then the ending term of the sequence (\ref{equ6.5}) is the simple module $V_i$. In this case, the same argument in the proof of Lemma \ref{lem4.3} shows that the sequence (\ref{equ6.5}) is almost split. We assume now that $2\leq j\leq n-1$. Since $M(\tau(i),j)\cong\textrm{DTr}M(i,j)$, by \cite[Proposition 2.2, ChV]{ARS}, we only need to verify that each non-isomorphism $f:M(i,j)\rightarrow M(i,j)$ factors through the following map: $$(id,0)\otimes id_{M(i,j)}:M(1,2)\otimes M(i,j)\rightarrow M(i,j).$$ Note that $M(i,j)$ is indecomposable and uniserial, and $f$ is not an isomorphism. Thus, the image of $f$ is contained in $JM(i,j)$, the unique maximal submodule of $M(i,j)$. Hence the left triangle in the following diagram is commutative:
$$
\xymatrix{
  M(i,j) \ar[rr]^{f} \ar[dr]_{f}
   & & M(i,j)     \\
   & JM(i,j) \ar@{^{(}->}[ur] \ar[rr]_{\iota}
   & & M(1,2)\otimes M(i,j) \ar[ul]_{(id,0)\otimes id_{M(i,j)}}.  }
$$
By Lemma \ref{lem3.3}, there exists an injective map $\iota$ from $JM(i,j)$ to $M(1,2)\otimes M(i,j)$ given by $\iota(x^kv)=1\otimes
x^kv-\zeta_kx\otimes
x^{k-1}v$ (see (\ref{equ3.2})), and
$$((id,0)\otimes id_{M(i,j)})(1\otimes x^{k}v-\zeta_{k}x\otimes x^{k-1}v)=1\otimes x^{k}v=x^{k}v,$$ for $1\leq k\leq j-1$ and $v\in V_i$.
This implies that the right triangle in the above diagram is also commutative. As a consequence, the non-isomorphism $f$ factors through the map $(id,0)\otimes id_{M(i,j)}$.\epf

\begin{prop}\label{prop3.2}For $1\leq i\leq m$, the following hold:
\begin{enumerate}
\item[(1)] $M(1,2)\otimes M(i,j)\cong M(i,j)\otimes M(1,2)\cong M(i,j+1)\oplus M(\tau(i),j-1)$, for $2\leq j\leq n-1$.
\item[(2)] $M(1,2)\otimes M(i,n)\cong M(i,n)\otimes M(1,2)\cong M(i,n)\oplus M(\tau(i),n)$.
\end{enumerate}
\end{prop}
\pf (1) For each non-projective indecomposable module $M(i,j)$, on the one hand, the sequences (\ref{equ6.5}) and (\ref{equ6.6}) ending with $M(i,j)$ are almost split. On the other hand, it follows from the proof of \cite[Theorem 2.1, ChVI]{ARS} that the following sequence is also almost split (we omit the maps of the sequence):
$$0\rightarrow M(\tau(i),j)\rightarrow M(\tau(i),j-1)\oplus M(i,j+1)\rightarrow M(i,j)\rightarrow0.$$
Since an almost split sequence is uniquely determined  by its beginning and ending terms, we obtain that $M(1,2)\otimes M(i,j)\cong M(i,j)\otimes M(1,2)\cong M(\tau(i),j-1)\oplus M(i,j+1)$.

Part (2) is obvious as the almost split sequences (\ref{equ6.5}) and (\ref{equ6.6}) with $j=n$ are split because $M(i,n)$ is projective. \epf

It is possible to give the decomposition of the tensor product $M(1,k)\otimes M(1,l)$ by virtue of Proposition \ref{prop3.2}. However, we do not continue this decomposition since it is more tedious. Instead, we leave it to the next section and express it as a multiplication rule in the Green ring of $H$.

To end this section, we describe the dual of the indecomposable
$H$-modules, which will be used later for the dual bases of the Green ring. Let $M$ be a finite dimensional $H$-module.
The dual space $M^*:=\textrm{Hom}(M,\mathbbm{k})$ is an $H$-module
given by $(hf)(v)=f(S(h)v),$ for $h\in H,\ f\in M^*$ and $v\in M.$

\begin{prop}\label{prop4.2}
For any $1\leq i\leq m$ and $1\leq j\leq n$,
\begin{enumerate}
\item[(1)] $M(i,j)^*\cong M(\tau^{1-j}(i^*),j)$, where $1\leq
i^*\leq m$ such that $(V_i)^*\cong V_{i^*}$.
\item[(2)] $M(i,j)^{**}\cong M(i,j)$.
\end{enumerate}
\end{prop}
\pf (1) We first study the dual $M(1,j)^*$ of the indecomposable
module $M(1,j)$. Let $\{1,x,\cdots,x^{j-1}\}$ be the basis of $M(1,j)$ and $\{1^*,x^*,\cdots,(x^{j-1})^*\}$ the dual basis of
$M(1,j)$. Then the actions of $h\in G$ and $y$ on the dual basis are given
respectively by $h(x^k)^*=\chi^{k}(h)(x^k)^*$ and
$y(x^k)^*=-q^{k-1}(x^{k-1})^*$, for $0\leq k\leq j-1$. Here
we point out that $y(1^*)=0$. Thus $M(1,j)^*$ is isomorphic to
$V_{\chi^{j-1}}\otimes M(1,j)$ and the isomorphism is given by
$$(x^k)^*\mapsto(-1)^kq^{\frac{k(k-1)}{2}}(w\otimes x^{j-1-k}),$$ for $0\leq k\leq j-1$, $0\neq w\in V_{\chi^{j-1}}$.
It follows from Proposition \ref{prop3.1} that
\begin{align*}
M(i,j)^*&\cong (M(1,j)\otimes V_i)^*\cong(V_i)^*\otimes
M(1,j)^*\\&\cong(V_i)^*\otimes V_{\chi^{j-1}}\otimes M(1,j)
\cong V_{i^*}\otimes V_{\chi^{j-1}}\otimes M(1,j)\\&\cong
V_{\tau^{1-j}(i^*)}\otimes M(1,j)\cong M(\tau^{1-j}(i^*),j).
\end{align*}

(2) Follows from  the fact that the square of the antipode of $H$ is inner.\epf

\section{\bf The structure of Green rings}

In this section, we present the Green ring $r(H_{\mathcal{D}})$ of a  pointed rank one Hopf algebra $H_{\mathcal{D}}$ with generators and relations. We show that the Green ring $r(H_{\mathcal{D}})$ is both Frobenius and symmetric.

Let $A$ be a Hopf algebra over a field $\kk$ and $F(A)$ the free abelian group generated by the isomorphism classes $[V]$ of finite dimensional
$A$-modules $V$. The abelian group $F(A)$ becomes a ring if we endow $F(A)$ with a multiplication given by the tensor product $[M][N]=[M\otimes N]$. The Green ring (or the representation ring) $r(A)$ of the Hopf algebra $A$ is defined to be the quotient ring of $F(A)$ modulo the relation $[M\oplus N]=[M]+[N]$. The identity of the associative ring  $r(A)$ is  represented by the trivial $A$-module $[\mathbbm{k}_\varepsilon]$. Note that
$r(A)$ has a $\mathbb{Z}$-basis consisting of the isomorphism classes of indecomposable modules, see \cite{DL,Wa,Wi}. The Grothendieck ring $G_0(A)$ of the Hopf algebra $A$ is the quotient ring of $F(A)$ modulo exact sequences   of $A$-modules $0\rightarrow M\rightarrow W\rightarrow N\rightarrow0$, i.e. $[W]=[M]+[N]$. The Grothendieck ring  $G_0(A)$ possesses a basis given by the isomorphism classes of simple modules. Both $r(A)$ and  $G_0(A)$ are  augmented $\mathbb{Z}$-algebras with the dimension augmentation. Moreover, there is  a natural ring epimorphism from $r(A)$ to $G_0(A)$. If $A$ is semisimple, then the ring epimorphism is the identity map.

Now let $r(H)$ be the Green ring of the Hopf algebra $H=H_{\mathcal {D}}$. We denote by $M[i,j]$ the isomorphism class of the indecomposable $H$-module $M(i,j)$ in $r(H)$. In particular, we set $1=[V_1]$ and
$a=[V_{\chi^{-1}}]=[V_{\tau(1)}]$.

\begin{prop}\label{prop4.1} The following hold in the Green ring $r(H)$:
\begin{enumerate}
\item[(1)] $M[i,j]=[V_i]M[1,j]=M[1,j][V_i]$, for $1\leq
i\leq m$ and $1\leq j\leq n$.
\item[(2)] $M[1,2]M[1,j]=M[1,j]M[1,2]=M[1,j+1]+aM[1,j-1]$,
for $2\leq j\leq n-1$.
\item[(3)] $M[1,2]M[1,n]=M[1,n]M[1,2]=(1+a)M[1,n]$.
\item[(4)] $r(H)$ is commutative and generated  by elements $[V_i]$, $1\leq i\leq m$ and $M[1,2]$ over $\mathbb{Z}$.
\end{enumerate}
\end{prop}
\pf Part (1) follows from Proposition \ref{prop3.1}. Part (2) and Part (3) follow from Proposition \ref{prop3.2}. Part (4) is a consequence of Part (1), Part (2) and the fact that $\{M[i,j]\mid 1\leq i\leq m,1\leq j\leq n\}$ forms a basis of $r(H)$.\epf

In the following, we give the decomposition rule of the product $M[1,k]M[1,l]$ in $r(H)$. This is correspondence to the decomposition of the tensor product $M(1,k)\otimes M(1,l)$ of $H$-modules.

\begin{prop}\label{col4.1}Let $1\leq k,l\leq n$. The following hold in $r(H)$:
\begin{enumerate}
\item[(1)] If $k+l-1\leq n$, then
$$M[1,k]M[1,l]=\sum_{t=0}^{\min\{k,l\}-1}M[\tau^t(1),k+l-1-2t].$$
\item[(2)] If $k+l-1\geq n$, let $r=k+l-1-n$. Then
$$M[1,k]M[1,l]=\sum_{t=0}^{r}M[\tau^t(1),n]+\sum_{t=r+1}^{\min\{k,l\}-1}M[\tau^t(1),k+l-1-2t].$$
\end{enumerate}
\end{prop}
\pf (1) We proceed by induction on $k+l-1$ for $1\leq k+l-1\leq n$.
It is obvious that the identity holds
for $k+l-1=1$. For a fixed $1<p\leq n-1$, suppose that the identity holds for $1<k+l-1\leq p$. We show that it holds for the case $k+l-1=p+1$. We may now assume that $k\geq 2$ without loss of generality. Since $k+l-1=p+1$ implies that $(k-1)+l-1\leq p$ and $(k-2)+l-1\leq p$, we may apply the induction hypothesis on $(k-1)+l-1\leq p$ and $(k-2)+l-1\leq p$. We obtain the following two equalities:
$$M[1,k-1]M[1,l]=\sum_{t=0}^{\min\{k-1,l\}-1}M[\tau^t(1),k-1+l-1-2t],$$
$$M[1,k-2]M[1,l]=\sum_{t=0}^{\min\{k-2,l\}-1}M[\tau^{t}(1),k-2+l-1-2t].$$
Now consider the product $ M[1,2](M[1,k-1]M[1,l])$. On the one hand,
we apply  the induction assumption to get:
\begin{align*} &\ \ \ \  M[1,2](M[1,k-1]M[1,l])\\
&= M[1,2]\sum_{t=0}^{\min\{k-1,l\}-1}M[\tau^t(1),k-1+l-1-2t]\\
&=\sum_{t=0}^{\min\{k-1,l\}-1}M[1,2]M[\tau^t(1),k-1+l-1-2t]\\
&=\sum_{t=0}^{\min\{k-1,l\}-1}(M[\tau^t(1),k-1+l-2t]+M[\tau^{t+1}(1),k-1+l-2-2t]).
\end{align*}
On the other hand, if we apply Proposition \ref{prop4.1} (2) on the product, we obtain:
\begin{align*}
&\ \ \ \  (M[1,2]M[1,k-1])M[1,l]\\
&=(M[1,k]+M[\tau(1),k-2])M[1,l]\\
&=M[1,k]M[1,l]+aM[1,k-2]M[1,l]\\
&=M[1,k]M[1,l]+\sum_{t=0}^{\min\{k-2,l\}-1}M[\tau^{t+1}(1),k-2+l-1-2t].
\end{align*}
These equations imply that \begin{align*}&\ \ \ \  M[1,k]M[1,l]+\sum_{t=0}^{\min\{k-2,l\}-1}M[\tau^{t+1}(1),k-2+l-1-2t]\\
&=\sum_{t=0}^{\min\{k-1,l\}-1}(M[\tau^t(1),k-1+l-2t]+M[\tau^{t+1}(1),k-1+l-2-2t]).\end{align*}
By discussing the cases $k-1<l$, $k-1=l$ and $k-1>l$, we  obtain that
$$M[1,k]M[1,l]=\sum_{t=0}^{\min\{k,l\}-1}M[\tau^t(1),k+l-1-2t].$$
Thus we have proved Part (1) for the case $k+l-1=p+1$.

(2) The proof is similar to the proof of Part (1) by induction on $k+l-1$ for $n\leq k+l-1\leq2n-1$. \epf

Now we are ready to give the structure of the Green ring $r(H)$. Let $\mathbbm{k}[y,z]$ be the polynomial ring with variables $y$ and
$z$ over $\mathbbm{k}$ and $F_i(y,z)$ the polynomials in $\mathbbm{k}[y,z]$ defined recursively as follows:
\be\label{equ7}F_1(y,z)=1,\ F_2(y,z)=z,\ F_i(y,z)=zF_{i-1}(y,z)-yF_{i-2}(y,z),\ i\geq3.\ee
Note that $F_i(y,z),\ i\geq1$ are called the generalized Fibonacci polynomial in \cite{CVZ, LZ} since if $y=-1$ then the polynomials $F_i(-1,z)\in\mathbbm{k}[z]$ are the well-known Fibonacci polynomials \cite{HB}. These polynomials are also referred to as the Dickson polynomials (of the second type), see \cite{LMT, HMSY}. As we shall see in the following, the Dickson polynomials are  fundamental factors in the structure of the Green ring $r(H)$.

\begin{theo}\label{th4.1}
Let $r(\mathbbm{k}G)$ be the Green ring of the group algebra $\kk G$ and $r(\mathbbm{k}G)[z]$ the polynomial ring with variable
$z$ over $r(\mathbbm{k}G)$. Then the Green ring $r(H)$ is isomorphic to $r(\mathbbm{k}G)[z]/I$, where $I$ is the
ideal of $r(\mathbbm{k}G)[z]$ generated by the element
$(1+a-z)F_n(a,z)$.
\end{theo}
\pf According to Proposition \ref{prop4.1}, $r(H)$ is generated as a ring by $M[1,2]$ over
$r(\mathbbm{k}G)$. Hence there is a unique ring epimorphism $\Phi$
from $r(\mathbbm{k}G)[z]$ to $r(H)$ such that
$$\Phi:r(\mathbbm{k}G)[z]\rightarrow r(H),\ g(z)\mapsto g(M[1,2]),$$ for any polynomial $g(z)\in r(\mathbbm{k}G)[z]$.
It is easy to check by induction on $j$ that
$\Phi(F_j(a,z))=M[1,j]$, for $1\leq j\leq n$. Now let $I$ be the
ideal of $r(\mathbbm{k}G)[z]$ generated by the element
$(1+a-z)F_n(a,z)$. By Proposition \ref{prop4.1}, we have
$$\Phi((1+a-z)F_n(a,z))=(1+a-M[1,2])M[1,n]=0.$$
This leads to a natural ring epimorphism $\overline{\Phi}$ from
$r(\mathbbm{k}G)[z]/I$ to $r(H)$ such that
$\overline{\Phi}(\overline{g(z)})=\Phi(g(z))$ for any $g(z)\in
r(\mathbbm{k}G)[z]$, where $\overline{g(z)}$ stands for the coset
$g(z)+I$ in $r(\mathbbm{k}G)[z]/I$. Observe that, as a
$\mathbb{Z}$-module, $r(\mathbbm{k}G)[z]/I$ has a $\mathbb{Z}$-basis
$\{\overline{[V_i]z^j}\mid1\leq i\leq m,\ 0\leq j\leq n-1\}$. Thus,
we conclude that $r(\mathbbm{k}G)[z]/I$ and $r(H)$
both have the same rank $mn$ as free $\mathbb{Z}$-modules and
hence the map $\overline{\Phi}$ is an isomorphism.\epf

\begin{rem} $(1)$ If $H$ is a Taft algebra or a generalized Taft algebra, then Theorem \ref{th4.1} recovers the main results of \cite{CVZ} and \cite{LZ}.

$(2)$ Since the Jacobson radical $J=(y)$ of $H$ is a Hopf ideal and $H/J\cong\mathbbm{k}G$, the Green ring
$r(\mathbbm{k}G)$ $($that is, the Grothendieck ring $G_0(\kk G))$ generated by $[V_i]$, $1\leq i\leq m$ is a subring of $r(H)$, and the Grothendieck ring $G_0(H)$ of $H$ is isomorphic
to $r(\mathbbm{k}G)$, see \cite{Lo}.

$(3)$ Consider $F_j(a,z)$, $1\leq j\leq n$ the Dickson polynomials in $r(\kk G)[z]$. Note that $\Phi(F_j(a,z))=M[1,j]$, see the proof of Theorem \ref{th4.1}.  By Proposition \ref{col4.1}, we have the following expressions for the multiplication of the Dickson polynomials of the second type.
\begin{itemize}
  \item If $k+l-1\leq n$, then
$$F_k(a,z)F_l(a,z)\equiv\sum_{t=0}^{\min\{k,l\}-1}a^tF_{k+l-1-2t}(a,z)\ (mod\ I).$$
  \item If $k+l-1\geq n$, let $r=k+l-1-n$. Then
$$F_k(a,z)F_l(a,z)\equiv\sum_{t=0}^{r}a^tF_n(a,z)+\sum_{t=r+1}^{\min\{k,l\}-1}a^tF_{k+l-1-2t}(a,z)\ (mod\ I),$$
\end{itemize}
where $I$ is the ideal given in Theorem \ref{th4.1}.
\end{rem}

For every finite dimensional pointed rank one Hopf algebra $H$ of nilpotent type, the computation of the Green ring of $H$ reduces to the computation of the Grothendieck ring of the group $G(H)$, where $G(H)$ is the group of the group-like elements of $H$. We now apply Theorem \ref{th4.1} to compute the Green ring of such a Hopf algebra $H$ with $G(H)$ being a dihedral group.

\begin{exa}
Assume that $s>0$ is a fixed odd integer. The dihedral group of
order $2s$ is defined as follows:
$$D_{2s}=\langle b,c\mid b^2=c^{2s}=(cb)^2=1\rangle,$$
where $c^s$ is the
unique non-trivial central element of $D_{2s}$. The simple modules
over the group algebra $\kk D_{2s}$ are given as follows (see e.g. \cite{Se}):
\begin{itemize}
  \item  Four simple modules of dimension 1:
$$F(i,j):c\mapsto(-1)^i,\ b\mapsto(-1)^j,\ i,j\in\mathbb{Z}_2,$$
  \item   $s-1$  simple modules of dimension 2:
$$V(l):c\mapsto\left(
                 \begin{array}{cc}
                   \theta^l & 0 \\
                   0 & \theta^{-l} \\
                 \end{array}
               \right)
,\ b\mapsto\left(
             \begin{array}{cc}
               0 & 1 \\
               1 & 0 \\
             \end{array}
           \right),\ 1\leq l\leq s-1,
$$ where  $\theta$ is a $2s$-th primitive root of unity.
\end{itemize}

The decompositions of the tensor products of simple modules
of $\kk D_{2s}$ are as follows (the computations are straightforward):
\begin{enumerate}
\item[$(1)$] $F(i,j)\otimes F(k,t)\cong F(i+k,j+t),\
i,j,k,t\in \mathbb{Z}_2$.
\item[$(2)$] $F(0,j)\otimes V(l)\cong
V(l),\ j\in \mathbb{Z}_2,\ 1\leq l\leq s-1$.
\item[$(3)$]
$F(1,j)\otimes V(l)\cong V(s-l),\ j\in \mathbb{Z}_2,\ 1\leq l\leq
s-1$.
\item[$(4)$] $V(l)\otimes V(h)\cong V(l-h)\oplus V(l+h),\
1\leq l,h\leq\frac{s-1}{2}\ \textrm{and}\ l>h$.
\item[$(5)$]
$V(l)\otimes V(l)\cong V(2l)\oplus F(0,0)\oplus F(0,1),\ 1\leq
l\leq\frac{s-1}{2}$.
\end{enumerate}
From the above decompositions, it follows that
$$[V(i)]=[V(i-1)][V(1)]-[V(i-2)],$$ for $2\leq i\leq s-1$, where $[V(0)]:=[F(0,0)]+[F(0,1)].$ Thus, the
Grothendieck ring $G_0(\kk D_{2s})$ of $\kk D_{2s}$ is
generated  by elements $[F(1,0)],[F(0,1)]$ and $[V(1)]$.
Hence, there is a unique ring epimorphism $\varphi$ from
$\mathbb{Z}[x_1,x_2,x_3]$ to $G_0(\kk D_{2s})$ such that
$$\varphi(x_1)=[F(1,0)],\ \varphi(x_2)=[F(0,1)]\ \textrm{and}\
\varphi(x_3)=[V(1)].$$ We define a sequence of polynomials
$f_i(x_2,x_3)$ in the polynomial ring $\mathbb{Z}[x_2,x_3]$ recursively as follows:
\begin{gather*}
f_0(x_2,x_3)=1+x_2,\ f_1(x_2,x_3)=x_3,\\
f_i(x_2,x_3)=x_3f_{i-1}(x_2,x_3)-f_{i-2}(x_2,x_3),\ \textrm{for}\
i\geq2.
\end{gather*}
Let $I_0$ be the ideal of $\mathbb{Z}[x_1,x_2,x_3]$ generated by the
polynomials
$x_1^2-1,x_2^2-1,x_2x_3-x_3,f_{\frac{s+1}{2}}(x_2,x_3)-x_1f_{\frac{s-1}{2}}(x_2,x_3)$. It is straightforward to verify that the Grothendieck ring $G_0(\kk D_{2s})$ is isomorphic to the quotient ring
$\mathbb{Z}[x_1,x_2,x_3]/I_0$.

Denote by $\chi$ the character of the simple module $F(1,0)$. It is obvious that the order of
$\chi$ is 2 and the order of $\chi(a^s)=(-1)^s=-1$ is also 2. Let $H_{\mathcal {D}}$ be the Hopf algebra stemming from
the group datum $\mathcal {D}=(D_{2s},\chi,c^s,0)$.
By Theorem \ref{th2.2}, the set
of indecomposable $H_{\mathcal {D}}$-modules consists of simple
modules and their projective covers. Let
$a=[V_{\chi^{-1}}]=[V_{\chi}]=[F(1,0)]$. Then, $a^2=1$. By Theorem
\ref{th4.1}, the Green ring $r(H_{\mathcal {D}})$ is isomorphic to $r(\kk D_{2s})[z]/I_1$, where $I_1$ is the ideal generated by $(1+a-z)z$. Thanks to the isomorphism $r(\kk  D_{2s})=G_0(\kk D_{2s})\cong\mathbb{Z}[x_1,x_2,x_3]/I_0$, the Green ring $r(H_{\mathcal {D}})$ is
isomorphic to the quotient ring $\mathbb{Z}[x_1,x_2,x_3,x_4]/I_2$,
where $I_2$ is the ideal of $\mathbb{Z}[x_1,x_2,x_3,x_4]$ generated by
the polynomials
$x_1^2-1,x_2^2-1,x_2x_3-x_3,f_{\frac{s+1}{2}}(x_2,x_3)-x_1f_{\frac{s-1}{2}}(x_2,x_3)$
and $x_4^2-x_1x_4-x_4$. \epf
\end{exa}

Two Hopf algebras are said to be {\it gauge equivalent} if their representation categories are tensor equivalent.
It is obvious that two gauge equivalent Hopf algebras possess the same Green ring. However, the converse is not true in general. In fact, the representation category of a Hopf algebra $A$ can not be solely determined by the
Green ring $r(A)$ of $A$, see e.g. \cite[Remark 1.8]{Wa}. In the following, we give a sufficient condition for two non-gauge equivalent Hopf algebras $H_{\mathcal{D}}$ to share the same Green ring.

\begin{prop}\label{prop4}
Let $\mathcal {D}=(G,\chi,g,0)$ and $\mathcal {D'}=(G,\chi',g,0)$ be
two group data of nilpotent type. Denote the set $G_k=\{h\in G\mid h^2=g^{-k}\}$ and the element $\vartheta=\sum_{k=0}^{n-1}\sum_{h\in G_{k}}(-1)^kh^{-k}g^{-\frac{k(1+k)}{2}}$ in $\mathbbm{k}G$. Suppose that $\chi(\vartheta)\neq\chi'(\vartheta)$
and there is an automorphism $\delta$ of $r(\mathbbm{k}G)$ such that
$\delta(a)=a'$, where $a=[V_{\chi^{-1}}]$ and $a'=[V_{\chi'^{-1}}]$.
Then the representation categories of Hopf algebras $H_{\mathcal {D}}$ and $H_{\mathcal {D'}}$ are not tensor equivalent. But the Green rings $r(H_{\mathcal {D}})$ and $r(H_{\mathcal {D'}})$ are isomorphic.
\end{prop}
\pf It is straightforward to check that $\chi(\vartheta)$ (resp. $\chi'(\vartheta)$) is the trace of the antipode of Hopf algebra $H_{\mathcal {D}}$ (resp. $H_{\mathcal {D'}}$). The condition $\chi(\vartheta)\neq\chi'(\vartheta)$ implies that the representation categories of Hopf algebras $H_{\mathcal {D}}$ and $H_{\mathcal {D'}}$ are not tensor equivalent since the trace of antipode of a Hopf algebra is a gauge invariant (i.e. invariant under tensor equivalence) (see \cite{KMN}).  By Theorem \ref{th4.1}, the Green ring
$r(H_{\mathcal {D}})$ (resp. $r(H_{\mathcal {D'}})$) is isomorphic
to $r(\mathbbm{k}G)[z]/I$, where $I$ is the ideal of
$r(\mathbbm{k}G)[z]$ generated by the element $(1+a-z)F_n(a,z)$
(resp. $(1+a'-z)F_n(a',z)$). Hence the automorphism $\delta$ of
$r(\mathbbm{k}G)$ such that $\delta(a)=a'$ induces an isomorphism
from $r(H_{\mathcal {D}})$ to $r(H_{\mathcal {D'}})$. \epf

\begin{exa}
Let $G$ be a cyclic group of order $n$ generated by $g$, $\chi$ a
linear character of $G$ such that the order of $\chi$ is $n$ and
$a=[V_{\chi^{-1}}]$. Given two group data $\mathcal
{D}=(G,\chi,g,0)$ and $\mathcal {D'}=(G,\chi^i,g,0)$ such that
$\textrm{gcd}(i,n)=1$. Then the Hopf algebras $H_{\mathcal {D}}$ and
$H_{\mathcal {D'}}$ are nothing but two Taft Hopf algebras. It is easy to see that the Green ring $r(\mathbbm{k}G)$ is isomorphic to
$\mathbbm{Z}C_n$, where $C_n=\langle a\rangle$ is a cyclic group of
order $n$. Let $\delta$ be the endomorphism of $\mathbbm{Z}C_n$
defined by $\delta(a)=a^i$. Then $\delta$ is an automorphism of
$\mathbbm{Z}C_n$ since $\textrm{gcd}(i,n)=1$. By Proposition \ref{prop4} (or \cite{Sch}) the representation categories of the two Taft Hopf algebras $H_\mathcal{D}$ and $H_\mathcal{D'}$ are not gauge equivalent since $\chi(g)\neq\chi^i(g)$. But the Green rings $r(H_\mathcal{D})$ and $r(H_\mathcal{D'})$ are
isomorphic.
\end{exa}

In the sequel, we show that the Green ring $r(H)$ possesses an associative symmetric and nonsingular $\mathbb{Z}$-bilinear form with a pair of dual bases associated to  almost split sequences.

Observe that the Green ring $r(H)$ is commutative and the dual modules given in Proposition \ref{prop4.2} induce an automorphism of $r(H)$ as follows: \be\label{equ4.4}M[i,j]^*=M[\tau^{1-j}(i^*),j]=a^{1-j}[V_i^*]M[1,j].\ee
This automorphism of $r(H)$ is an involution. Namely, $M[i,j]^{**}=M[i,j]$, for $1\leq i\leq m$, $1\leq j\leq n$.

For the short exact sequence  (\ref{equ6.5}) ending with the non-projective indecomposable module $M(i,j)$, we follow the notations given in \cite[Section 4, ChVI]{ARS} and denote by $\delta_{M[i,j]}$, $1\leq i\leq m$, $1\leq j\leq n-1$ the following elements in $r(H)$:
\begin{align*}\delta_{M[i,j]}
&=M[\tau(i),j]-[M(1,2)\otimes M(i,j)]+M[i,j]\\
&=aM[i,j]-M[1,2]M[i,j]+M[i,j]\\
&=(1+a-M[1,2])M[i,j].
\end{align*}
For the indecomposable projective module $M(i,n)$, $1\leq i\leq m$, we write
\begin{align*}\delta_{M[i,n]}&=M[i,n]-[JM(i,n)]\\
&=M[i,n]-M[\tau(i),n-1]\\
&=M[i,n]-aM[i,n-1].
\end{align*}
We define  $$\langle[V],[W]\rangle=:\dim \textrm{Hom}_H(V,W),$$ for any two $H$-modules $V$ and $W$. Then $\langle\ ,\ \rangle$ induces a $\mathbb{Z}$-bilinear form over $r(H)$. Note that $M(i,j)$ is finite dimensional indecomposable and the endomorphism ring of $M(i,j)$ is local. Thus, $\textrm{End}_H(M(i,j))/\textrm{rad End}_H(M(i,j))$ is of dimension 1,
as the field $\kk$ is assumed to be algebraically closed. The following result can be found from \cite[Proposition 4.1, ChVI]{ARS}.

\begin{lem}\label{prop4.3}For any indecomposable $H$-module $V$ and $x\in r(H)$, the following hold:
\begin{enumerate}
\item[(1)] $\langle[V],\delta_{M[i,j]}\rangle=1$ if and only if $V\cong M(i,j)$.
\item[(2)] $x=\sum_{i=1}^m\sum_{j=1}^n\langle x,\delta_{M[i,j]}\rangle M[i,j]$.
\item[(3)] $\{\delta_{M[i,j]}\mid 1\leq i\leq m,\ 1\leq j\leq n\}$ forms a basis of $r(H)$.
\end{enumerate}
\end{lem}

The $\mathbb{Z}$-bilinear form $\langle\ ,\ \rangle$ defined above is neither associative nor symmetric. However, we may modify it as follows:
$$([V],[W])=\langle[V],[W]^*\rangle=\dim \textrm{Hom}_H(V,W^*),$$ for any $H$-modules $V$ and $W$. Then $(\ ,\ )$ extends to a $\mathbb{Z}$-bilinear form over $r(H)$. It follows from the canonical $H$-module isomorphisms: $\textrm{Hom}_H(V,W)\cong (W\otimes V^*)^H$, $V\cong V^{**}$ and the $\kk$-linear isomorphism: $\textrm{Hom}_H(V^*,W^*)\cong \textrm{Hom}_H(W,V)$,  that the $\mathbb{Z}$-bilinear form $(\ ,\ )$ is associative and symmetric.
Denote by  $\delta_{M[i,j]}^*$ the image of $\delta_{M[i,j]}$ under the dual automorphism.
Then by (\ref{equ4.4}), we have
$$\delta_{M[i,j]}^*=\delta_{M[\tau^{-j}(i^*),j]},\ \textrm{for}\ 1\leq i\leq m,\ 1\leq j\leq n-1,$$
$$\delta_{M[i,n]}^*=a^{1-n}[V_i^*](M[1,n]-M[1,n-1]),\ \textrm{for}\ 1\leq i\leq m.$$
Now the Lemma \ref{prop4.3} can be modified as follows.

\begin{lem}\label{prop4.5}For any indecomposable $H$-module $V$ and $x\in r(H)$, we have the following:
\begin{enumerate}
\item[(1)] $([V],\delta_{M[i,j]}^*)=1$ if and only if $V\cong M(i,j)$.
\item[(2)] $x=\sum_{i=1}^m\sum_{j=1}^n(x,\delta_{M[i,j]}^*)M[i,j]$.
\item[(3)] $\{\delta_{M[i,j]}^*\mid 1\leq i\leq m,\ 1\leq j\leq n\}$ is a basis of $r(H)$.
\end{enumerate}
\end{lem}

Note that $M[i,j]=[V_i](1+a+\cdots+a^{j-1})$ holds in the Grothendieck ring $G_0(H)$ of $H$, for $1\leq i\leq m$ and $1\leq j\leq n$. In view of this, the natural ring epimorphism from $r(H)$ to $G_0(H)$ is given by $$M[i,j]\mapsto[V_i](1+a+\cdots+a^{j-1}),\ \textrm{for}\ 1\leq i\leq m,\ 1\leq j\leq n,$$
whose kernel is precisely spanned by the basis elements: $\delta_{M[i,j]}$ (or $\delta_{M[i,j]}^*$), for $1\leq i\leq m$ and $1\leq j\leq n-1$, see \cite[Theorem 4.4, ChVI]{ARS}.

For each simple $H$-module $V_i$, we denote by $\chi_i$ the character of the simple $\kk G$-module $V_i$. It is well known that the free abelian group $\textrm{Irr}(\kk G):=\sum_{i=1}^m\mathbb{Z}\chi_i$ with the convolution product forms a ring, called {\it the character ring} of $\kk G$. The character ring $\textrm{Irr}(\kk G)$ is isomorphic to the Green ring $r(\kk G)$ via the map $\chi_i\mapsto [V_i]$, for $1\leq i\leq m$. The equality
\be\label{equ4.5}[V_i][V_j]=\delta_{i,j^*}[V_1]+\sum_{k=2}^{m}\gamma_k[V_k]\ee follows from the fact that
$\chi_i\chi_j=\delta_{i,j^*}\chi_1+\sum_{k=2}^{m}\gamma_k\chi_k$, for some non-negative integers $\gamma_k$.
The equality (\ref{equ4.5}) could also be interpreted by Lemma \ref{prop4.5}. That is
\begin{align*}[V_i][V_j]&=\sum_{k=1}^m\sum_{l=1}^n([V_i][V_j],\delta_{M[k,l]}^*)M[k,l]
=\sum_{k=1}^m([V_i][V_j],\delta_{M[k,1]}^*)M[k,1]\\
&=([V_i][V_j],\delta_{M[1,1]}^*)M[1,1]+\sum_{k=2}^m([V_i][V_j],\delta_{M[k,1]}^*)M[k,1]\\
&=([V_i],\delta_{M[j^*,1]}^*)M[1,1]+\sum_{k=2}^m([V_i][V_j],\delta_{M[k,1]}^*)M[k,1]\\
&=\delta_{i,j^*}[V_1]+\sum_{k=2}^m([V_i][V_j],\delta_{M[k,1]}^*)[V_k].
\end{align*}

It follows from Lemma \ref{prop4.5} (2) that the bilinear form $(\ ,\ )$ is nonsingular. Thus the Green ring $r(H)$ is symmetric as a $\mathbb{Z}$-algebra. As a consequence, various ring theoretic properties of Frobenius and symmetric algebras always hold for $r(H)$. For instance, the Green ring $r(H)$ has the following properties, see \cite{KS, Ka, Lor} for related materials.

\begin{prop} The following hold in $r(H)$.
\begin{enumerate}
\item The bases $\{M[i,j]\mid1\leq i\leq m, 1\leq j\leq n\}$ and $\{\delta_{M[i,j]}^*\mid 1\leq i\leq m, 1\leq j\leq n\}$, form a dual pair of bases for $r(H)$ with respt to the $\mathbb{Z}$-bilinear form $(\ ,\ )$.
\item  The Frobenius homomorphism $\phi:r(H)\rightarrow \mathbb{Z}$ is given by $\phi(x)=(x,1)$, for any $x\in r(H)$.
\item The Casimir element of $r(H)$ is given by $\sum_{i=1}^m\sum_{j=1}^n\delta_{M[i,j]}^*\otimes M[i,j]=\sum_{i=1}^m\sum_{j=1}^n M[i,j]\otimes\delta_{M[i,j]}^*$.
\item The integral $\int_{r(H)}=\{b\in r(H)\mid bx=\dim(x)b,\ \textrm{for}\ x\in r(H)\}$ determined by the dimension augmentation is an ideal of $r(H)$ generated by $[H]$, where $[H]$ is the isomorphism class of regular representation of $H$. Moreover, the rank of $\int_{r(H)}$ is 1 and $([H],\ )$=dim.
\end{enumerate}
\end{prop}

\section{\bf Jacobson radicals of Green rings}
In this section, we study the Jacobson radical $J(r(H))$ of the Green ring $r(H)$. It turns out that the Jacobson radical of $r(H)$ is a principal ideal no matter what the group $G$ is. We assume that the Hopf algebra $H=H_\mathcal {D}$ is defined over $\mathbbm{k}=\mathbb{C}$.


Suppose that the order of $\chi$ in the group datum $\mathcal
{D}=(G,\chi,g,0)$ is $l$. Then so is the order of
$a=[V_{\chi^{-1}}]$. Note that $q=\chi(g)$ is a primitive $n$-th
root of unity, and $q^l=(\chi(g))^l=\chi^l(g)=1$. Hence, $l$ is
divisible by $n$. Since $\mathbb{C}G$ is semisimple, the
complexified Green ring
$R(\mathbb{C}G):=\mathbb{C}\otimes_{\mathbb{Z}}r(\mathbb{C}G)$ is a commutative semisimple algebra, see
\cite{Wi,Zhu}. As a consequence, all the simple modules over
$R(\mathbb{C}G)$ is one dimensional and the number of
non-isomorphic simple modules is equal to the rank of
$r(\mathbb{C}G)$, which is equal to
the number of non-isomorphic simple $\mathbb{C}G$-modules, namely,
$m$. Let $W_1,W_2,\cdots,W_m$ be a complete set of non-isomorphic
simple $R(\mathbb{C}G)$-modules. Then the action of $a$ on
$W_j$ is a scalar multiplication by some $\omega^{t_j}$, where
$\omega=\cos\frac{2\pi}{l}+i\sin\frac{2\pi}{l}$ is a primitive
$l$-th root of unity and $0\leq t_j\leq l-1$.
We divide the index set into three parts:
$$\Omega_1=\{j\mid 1\leq j\leq m, t_j=0\},$$
$$\Omega_2=\{j\mid 1\leq j\leq m, t_j\neq0\ \textrm{and}\ \frac{l}{n}\nmid
t_j\}$$ and $$\Omega_3=\{j\mid 1\leq j\leq m, t_j\neq0\
\textrm{and}\ \frac{l}{n}\mid t_j\}.$$ The cardinality of the above
sets $\Omega_1$, $\Omega_2$ and $\Omega_3$ are denoted  $d_1$,
$d_2$ and $d_3$ respectively. Obviously, $d_1+d_2+d_3=m$. For $1\leq
j\leq m$, let $N_j$ stand for the number of distinct roots of the
equation $(z-\omega^{t_j}-1)F_n(\omega^{t_j},z)=0$. We then have the
following lemma, similar to \cite[Proposition 4.2]{LZ}.

\begin{lem}\label{lem4.2}
If $j\in\Omega_3$, then $N_j=n-1$. If $j\in \Omega_1\cup\Omega_2$, then $N_j=n$.
\end{lem}
\pf Let
$b=\cos(\frac{t_j\pi}{l}+\frac{3\pi}{2})+i\sin(\frac{t_j\pi}{l}+\frac{3\pi}{2})$. Then $b^2=-\omega^{t_j}$. The connection between the polynomials
$F_s(\omega^{t_j},z)$ and the Fibonacci polynomials $F_s(-1,z)$ is
established by induction on $s$ as follows:
$$F_s(\omega^{t_j},z)=b^{s-1}F_s(-1,b^{-1}z),\ \textrm{for}\ s\geq1.$$
In particular, $F_n(\omega^{t_j},z)=b^{n-1}F_n(-1,b^{-1}z)$. Note that
the distinct roots of the equation $F_n(-1,z)=0$ are $z_k=2i\cos\frac{k\pi}{n}$,
for $1\leq k\leq n-1$, here $i^2=-1$, see \cite{HB}. It follows that
the distinct roots of $F_n(\omega^{t_j},z)=0$ are
$z_k=2bi\cos\frac{k\pi}{n}$, for $1\leq k\leq n-1$. This implies
that the equation $(z-\omega^{t_j}-1)F_n(\omega^{t_j},z)=0$ has
roots $\omega^{t_j}+1$ and $2bi\cos\frac{k\pi}{n}$, $1\leq k\leq
n-1$. Now $\omega^{t_j}+1=2bi\cos\frac{k\pi}{n}$ if and only if
$\cos\frac{t_j\pi}{l}=\cos\frac{k\pi}{n}$ if and only if $k=s$ and
$t_j=\frac{ls}{n}$, for a unique $1\leq s\leq n-1$. Therefore, there
are $n-1$ distinct roots if $j\in\Omega_3$, and $n$, otherwise.\epf

\begin{prop}\label{prop4.4} Let $R(H):=\mathbb{C}\otimes_{\mathbb{Z}}r(H)$ be the
complexified Green ring of $r(H)$.
\begin{enumerate}
\item[$(1)$] $R(H)$ has exactly $mn-d_3$ simple modules and each
of them is of dimension one;
\item[$(2)$] The dimension of the Jacobson radical $J(R(H))$ of $R(H)$ is $d_3$.
\end{enumerate}
\end{prop}
\pf (1) The fact that $R(H)$ is commutative and the ground field is
$\mathbb{C}$ implies that each simple $R(H)$-module is of
dimension one. According to Theorem \ref{th4.1}, the Green ring $r(H)$
of $H$ is isomorphic to $r(\mathbb{C}G)[z]/I$. This means that each
one dimensional $R(H)$-module is a lift from a one dimensional
$R(\mathbb{C}G)[z]$-module. Since the action of $a$ on each simple
$R(\mathbb{C}G)$-module $W_j$ is a scalar multiplication by some
$\omega^{t_j}$, the $R(\mathbb{C} G)$-module $W_j$ becomes a simple $R(H)$-module if and only if the action of $z$ on $W_j$ is a scalar multiplication by a
root of the equation $(z-\omega^{t_j}-1)F_n(\omega^{t_j},z)=0.$ This
equation has $N_j$ distinct roots. We can conclude that the number of
non-isomorphic simple $R(H)$-modules lifted from $W_j$ is $N_j$. By
Lemma \ref{lem4.2}, the number of non-isomorphic simple
$R(H)$-modules is $\Sigma_{j=1}^mN_j=(m-d_3)n+d_3(n-1)=mn-d_3$.

(2) Note that the dimension of $R(H)$ is $mn$ and the dimension of
the quotient algebra $R(H)/J(R(H))$ is equal to the number of
non-isomorphic simple $R(H)$-modules, namely, $mn-d_3$. Thus, the
dimension of the Jacobson radical $J(R(H))$ of $R(H)$ is $d_3$. \epf

Now let $\theta$ be the following element in $r(\mathbb{C} G)$:
$$\theta=(1-a)(1+a^n+a^{2n}+\cdots+a^{(\frac{l}{n}-1)n}),$$
and $\langle\theta\rangle$ the principal ideal of $r(\mathbb{C}G)$ generated by $\theta$. Then $\mathbb{C}\otimes_{\mathbb{Z}}\langle\theta\rangle$
is an ideal of $R(\mathbb{C}G)=\mathbb{C}\otimes_{\mathbb{Z}}r(\mathbb{C}G)$.
Multiplying  $\langle \theta\rangle$ by the element $M[1,n]$, we get a $\mathbb{Z}$-submodule $M[1,n]\langle\theta\rangle$  of $r(H)$. Since $r(H)$ is free over $\mathbb{Z}$, the submodule $M[1,n]\langle\theta\rangle$ is  free as well.

\begin{lem}\label{lem4.4} The following holds:
\begin{enumerate}
\item[$(1)$] The rank of $\langle\theta\rangle$ is $d_3$.
\item[$(2)$] The rank of
$M[1,n]\langle\theta\rangle$ is equal to the
rank of $\langle\theta\rangle$.
\end{enumerate}
\end{lem}
\pf (1) Note that the quotient algebra
$R(\mathbb{C}G)/(\mathbb{C}\otimes_{\mathbb{Z}}\langle\theta\rangle)$
is commutative semisimple. We first determine the
dimension of
$R(\mathbb{C}G)/(\mathbb{C}\otimes_{\mathbb{Z}}\langle\theta\rangle)$
by calculating the number of non-isomorphic simple
$R(\mathbb{C}G)/(\mathbb{C}\otimes_{\mathbb{Z}}\langle\theta\rangle)$-modules.
Observe that each simple
$R(\mathbb{C}G)/(\mathbb{C}\otimes_{\mathbb{Z}}\langle\theta\rangle)$-module
is precisely a simple $R(\mathbb{C}G)$-module $W_j$ such that
$\theta W_j=0$, whereas, $\theta W_j=0$ if and only if
$$(1-\omega^{t_j})(1+\omega^{nt_j}+\omega^{2nt_j}\cdots+\omega^{(\frac{l}{n}-1)nt_j})=0,$$
if and only if  $j\in\Omega_1\cup\Omega_2$. As a consequence,
there exist exactly $d_1+d_2$ distinct simple
$R(\mathbb{C}G)/(\mathbb{C}\otimes_{\mathbb{Z}}\langle\theta\rangle)$-modules
and hence the dimension of
$R(\mathbb{C}G)/(\mathbb{C}\otimes_{\mathbb{Z}}\langle\theta\rangle)$
is $d_1+d_2$. It follows  that the dimension of
$\mathbb{C}\otimes_{\mathbb{Z}}\langle\theta\rangle$ is
$m-(d_1+d_2)=d_3$. Therefore, as a free $\mathbb{Z}$-module, the
rank of $\langle\theta\rangle$ is $d_3$.

(2) We prove the general case: if $I$ is an ideal of
$r(\mathbb{C}G)$, then both $I$ and $M[1,n]I$ have the same rank as free $\mathbb{Z}$-modules. Suppose that $\theta_1,\theta_2,\cdots,\theta_k$ form
a $\mathbb{Z}$-basis of $I$. It is obvious that $M[1,n]I$ is generated as a
$\mathbb{Z}$-module by $M[1,n]\theta_1,M[1,n]\theta_2,\cdots,M[1,n]\theta_k$.
We claim that the foregoing generators form a $\mathbb{Z}$-basis of
$M[1,n]I$. Indeed, if $\alpha_1M[1,n]\theta_1+\alpha_2M[1,n]\theta_2+\cdots+\alpha_kM[1,n]\theta_k=0$, where each $\alpha_i\in\mathbb{Z}$, then we have
\be\label{equ4.1} M[1,n](\alpha_1\theta_1+\alpha_2\theta_2+\cdots+\alpha_k\theta_k)=0.\ee
Denote by $K_0(H)$ the Grothendieck group of the category of finite
dimensional projective left $H$-modules, that is, the abelian group
generated by the isomorphism classes $M[i,n]$ of projective
$H$-modules $M(i,n)$ modulo the relations $[M(i,n)\oplus
M(j,n)]=M[i,n]+M[j,n]$. Then $K_0(H)$ is a free abelian group with
a basis $\{M[i,n]\mid1\leq i\leq m\}$. It is obvious that $K_0(H)$ admits a right action from $G_0(H)=r(\mathbb{C}G)$: $M[i,n]\cdot [V_j]=[M(i,n)\otimes V_j]$. Thus $K_0(H)$ is a right module over
$G_0(H)$. In fact, $K_0(H)$ is free of rank 1 with generator
$M[1,n]$ as a right $G_0(H)$-module. Now the equation
(\ref{equ4.1}) implies that
$$K_0(H)(\alpha_1\theta_1+\alpha_2\theta_2+\cdots+\alpha_k\theta_k)=0.$$
Hence, $\alpha_1\theta_1+\alpha_2\theta_2+\cdots+\alpha_k\theta_k=0$ since the right $G_0(H)$-module $K_0(H)$ is faithful, see \cite{Lo}. We obtain that $\alpha_i=0$, for $1\leq i\leq k$. Therefore, the rank
of $I$ is equal to the rank of $M[1,n]I$. \epf

\begin{theo}\label{th4.5}
The Jacobson radical of $r(H)$ is a principal ideal generated by the element $M[1,n]\theta$.
\end{theo}
\pf Note that $(1+a+\cdots+a^{n-1})\theta=(1-a^n)(1+a^n+a^{2n}\cdots+a^{(\frac{l}{n}-1)n})=1-a^l=0,$ and by Proposition \ref{col4.1}, $M[1,n]^2=(1+a+\cdots+a^{n-1})M[1,n].$
This yields that $(M[1,n]\theta)^2=0$. Hence
$M[1,n]\langle \theta\rangle\subseteq J(r(H))$.
Now by Lemma \ref{lem4.4} we know that the rank of
$M[1,n]\langle\theta\rangle$ is $d_3$, and
consequently, the rank of $J(r(H))$ is equal or greater than $d_3$.
On the other hand, since
$\mathbb{C}\otimes_{\mathbb{Z}}J(r(H))\subseteq J(R(H))$ and by
Proposition \ref{prop4.4}, the dimension of $J(R(H))$ is $d_3$. It
follows that the rank of $J(r(H))$ is equal or less than $d_3$. Thus, we
conclude that the rank of $J(r(H))$ is equal to $d_3$. Therefore,
$J(r(H))=M[1,n]\langle\theta\rangle$. Let $\langle M[1,n]\theta\rangle$ be the principal ideal of $r(H)$ generated by $M[1,n]\theta$. We have the following inclusions:
$J(r(H))=M[1,n]\langle \theta\rangle\subseteq\langle M[1,n]\theta\rangle\subseteq J(r(H)).$
It follows that $J(r(H))=\langle M[1,n]\theta\rangle$, as desired. \epf

\begin{rem}
In case $H$ is a generalized Taft algebra, the Jacobson radical of the Green ring $r(H)$ has been calculated in \cite{LZ}, where each nilpotent element is represented by a linear combination of certain projective indecomposables. Now we understand the form.
\end{rem}

\section{\bf Green rings of stable categories}

In this Section, we study the Green ring of the stable (module) category of a pointed rank one Hopf algebra $H$ of nilpotent type. We show that the complexified stable Green algebra is a group-like algebra and hence a bi-Frobenius algebra.

Let $H$-mod be the category of (finite dimensional) $H$-modules. Recall that the stable category $H$-\underline{mod} is the quotient category  of $H$-mod modulo the morphism factoring through the projective modules. This category is triangulated \cite{Ha} with the monoidal structure derived from that of $H$-mod. The Green ring of the stable category $H$-\underline{mod} is called  {\it the stable Green ring of $H$}, denoted $r_{st}(H)$.

\begin{prop}\label{prop6.1}
The stable Green ring $r_{st}(H)$  is isomorphic to the quotient ring $r(H)/\mathcal{P}$, where $\mathcal{P}$ is the ideal of $r(H)$ generated by the isomorphism classes of indecomposable projective $H$-modules. Precisely, $r_{st}(H)$ is isomorphic to $r(\kk G)[z]/(F_n(a,z))$.
\end{prop}
\pf The functor $F$ from $H$-mod to $H$-\underline{mod} given by $F(M)=M$ for any $H$-module $M$, and $F(\phi)=\underline{\phi}$ for $\phi\in\textrm{Hom}_{H}(M,N)$ with the canonical image $\underline{\phi}\in\textrm{\underline{Hom}}(M,N)$ is a full, dense tensor functor. Such a functor defines a ring epimorphism $f$ from $r(H)$ to $r_{st}(H)$ such that $f(\mathcal{P})=0$. Hence there is a unique ring epimorphism $\overline{f}$ from $r(H)/\mathcal{P}$ to $r_{st}(H)$ such that $\overline{f}(\overline{x})=f(x)$, for any $x\in r(H)$ with the the canonical image $\overline{x}\in r(H)/\mathcal{P}$. The rank of $r_{st}(H)$ is the same as the rank of $r(H)/\mathcal{P}$ since there is one to one correspondence between the indecomposable objects in $H$-\underline{mod} and the non-projective indecomposable modules in $H$-mod. We conclude that $r_{st}(H)$ is isomorphic to $r(H)/\mathcal{P}$. Since every projective indecomposable module $M(i,n)$ is isomorphic to $V_i\otimes M(1,n)$ (see Prop. 4.1), the ideal  $\mathcal{P}$ is a principal ideal generated by $M[1,n]$, i.e. $F_n(a,z)$. It follows that the stable Green ring is isomorphic to $r(\kk G)[z]/(F_n(a,z))$.
\epf

The stable Green ring $r_{st}(H)$ is semisimple (see Theorem \ref{th4.5}) with a $\mathbb{Z}$-basis $\{\overline{M[i,j]}\mid1\leq i\leq m,1\leq j\leq n-1\}$. Consider the complexified stable Green ring $R_{st}(H):=\mathbb{C}\otimes_{\mathbb{Z}}r_{st}(H)$, a semisimple algebra. We show that this algebra is both a group-like algebra and a bi-Frobenius algebra.  The notion of a group-like algebra was introduced by Doi in \cite{Doi1} generalizing the group algebra of a finite group and a scheme ring (Bose-Mesner algebra) of a non-commutative association scheme.

\begin{defi} Let $(A,\varepsilon,\mathbf{B},*)$ be a quadruple, where $A$ is a finite dimensional algebra over a field $\kk$ with unit 1, $\varepsilon$ is an algebra morphism from $A$ to $\kk$, the set $\mathbf{B}=\{b_i\}_{i\in I}$ is a $\kk$-basis of $A$ such that $0\in I$ and $b_0=1$, and $*$ is an involution of the index set $I$. Then $(A,\varepsilon,\mathbf{B},*)$ is called a group-like algebra if the following hold:
\begin{enumerate}
\item[(G1)] $\varepsilon(b_i)=\varepsilon(b_{i^{*}})\neq0$ for all $i\in I$.
\item[(G2)] $p_{ij}^k=p_{j^*i^*}^{k^*}$ for all $i,j,k\in I$, where $p_{ij}^k$ is the structure constant for $\mathbf{B}$ defined by $b_ib_j=\sum_{k\in I}p_{ij}^kb_k$.
\item[(G3)] $p_{ij}^0=\delta_{i,j^*}\varepsilon(b_i)$ for all $i,j\in I$.
\end{enumerate}
\end{defi}

The complexified Green ring $R(\kk G)=\mathbb{C}\otimes_{\mathbb{Z}}r(\kk G)$ is a group-like algebra. The algebra morphism $\varepsilon$ from $R(\kk G)$ to $\mathbb{C}$ is given by $\varepsilon([V_i])=\dim(V_i),$ for $1\leq i\leq m$. The basis set $\mathbf{B}$ can be chosen as the set: $$\{\dim(V_1)[V_1],\dim(V_2)[V_2],\cdots,\dim(V_m)[V_m]\},$$
 and the involution $*$ is induced by the dual map $[V_{i^*}]=[V_i]^*$, for $1\leq i\leq m$.

Let $\varepsilon$ be the algebra morphism from $R(\kk G)[z]$ to $\mathbb{C}$ such that $$\varepsilon([V_i])=\dim(V_i),\ \textrm{for}\ 1\leq i\leq m,\ \varepsilon(z)=2\cos\frac{\pi}{n}.$$ It follows from the proof of Lemma \ref{lem4.2} that $\varepsilon(F_j(a,z))=F_j(1,2\cos\frac{\pi}{n})\neq0$ for $1\leq j\leq n-1$ and $\varepsilon(F_n(a,z))=F_n(1,2\cos\frac{\pi}{n})=0$. Since $R_{st}(H)$ is isomorphic to $R(\kk G)[z]/(F_n(a,z))$ (see Proposition \ref{prop6.1}), $\varepsilon$ is exactly an algebra morphism from $R_{st}(H)$ to $\mathbb{C}$ such that $$\varepsilon(\overline{M[i,j]})=\dim(V_i)F_j(1,2\cos\frac{\pi}{n}).$$
Write $I=\{(i,j)\mid1\leq i\leq m,\ 1\leq j\leq n-1\}$ and $$b_{(i,j)}=\varepsilon(\overline{M[i,j]})\overline{M[i,j]},\ \textrm{for}\ (i,j)\in I.$$ Then the set $\mathbf{B}=\{b_{(i,j)}\mid (i,j)\in I\}$ forms a basis of $R_{st}(H)$ with $b_{(1,1)}=1$. Since the automorphism $*$ of $r(H)$ given by (\ref{equ4.4}) preserves the ideal $\mathcal{P}$, it induces an automorphism over the algebra $R_{st}(H)$. Namely, for each $(i,j)\in I$, \be\label{equ5.2}(\overline{M[i,j]})^*=\overline{M[i,j]^*}.\ee
Moreover, by (\ref{equ4.4}),
\begin{align}\label{equ5.1}\varepsilon(\overline{M[i,j]^*})&=\varepsilon(\overline{M[\tau^{1-j}(i^*),j]})
=\varepsilon(\overline{a^{1-j}})\varepsilon(\overline{[V_{i^*}]})\varepsilon(\overline{M[1,j]})\nonumber\\
&=\dim(V_{i^*})\varepsilon(\overline{M[1,j]})=\dim(V_i)\varepsilon(\overline{M[1,j]})\\
&=\varepsilon(\overline{M[i,j]}).\nonumber
\end{align}
Also, the the automorphism $*$ of $r(H)$ induces an involution on $I$. Namely, \be\label{equ5.3}(i,j)^*=(\tau^{1-j}(i^*),j),\ee for any $(i,j)\in I$. With the above notations, we have the following:

\begin{prop}
The quadruple $(R_{st}(H),\varepsilon,\mathbf{B},*)$ is a group-like algebra.
\end{prop}
\pf The condition (G1) follows from (\ref{equ5.1}) and (\ref{equ5.3}).

 For the condition (G2), we have  \be\label{equ5.4}(b_{(i,j)})^*=\varepsilon(\overline{M[i,j]})(\overline{M[i,j]})^*=\varepsilon(\overline{M[i,j]^*})\overline{M[i,j]^*}=b_{(i,j)^*},\ee
for any $(i,j)\in I$, because of (\ref{equ5.1}), (\ref{equ5.2}) and (\ref{equ5.3}). Now for any $(i,k),(j,l)\in I$, suppose that $b_{(i,k)}b_{(j,l)}=\sum_{(s,t)\in I}p_{(i,k)(j,l)}^{(s,t)}b_{(s,t)}$, where  $p_{(i,k)(j,l)}^{(s,t)}\in\mathbb{C}$. Applying the dual automorphism $*$ and (\ref{equ5.4}) to the forgoing equality, we obtain $b_{(j,l)^*}b_{(i,k)^*}=\sum_{(s,t)\in I}p_{(i,k)(j,l)}^{(s,t)}b_{(s,t)^*}$. On the other hand, $b_{(j,l)^*}b_{(i,k)^*}=\sum_{(r,q)\in I}p_{(j,l)^*(i,k)^*}^{(r,q)}b_{(r,q)}$. We conclude that $p_{(i,k)(j,l)}^{(s,t)}=p_{(j,l)^*(i,k)^*}^{(s,t)^*}$ for any $(i,k),(j,l),(s,t)\in I$.

Now we verify the condition (G3). For any $(i,k),(j,l)\in I$ such that $k\neq l$, note that $(j,l)^*=(\tau^{1-l}(j^*),l)$. By Proposition \ref{col4.1}, the multiplication $b_{(i,k)}b_{(j,l)}$ does not contain the term $b_{(1,1)}$. In this case, $p_{(i,k)(j,l)}^{(1,1)}=0=\delta_{(i,k),(j,l)^*}\varepsilon(b_{(i,k)})$ and the condition (G3) holds. If $k=l$, then by Proposition \ref{col4.1}, we have
\begin{align*}b_{(i,k)}b_{(j,k)}&=\varepsilon(\overline{M[i,k]})\varepsilon(\overline{M[j,k]})\overline{M[i,k]}\overline{M[j,k]}\\
&=\varepsilon(\overline{M[i,k]})\varepsilon(\overline{M[j,k]})\overline{[V_i][V_j]M[1,k]^2}\\
&=\varepsilon(\overline{M[i,k]})\varepsilon(\overline{M[j,k]})\overline{[V_i][V_j]}(\overline{a}^{k-1}+\sum_{s=2}^{k-1}\mu_s\overline{M[1,s]}),
\end{align*} for some non-negative integers $\mu_s$.
It follows from Equation (\ref{equ4.5}) that
$$\overline{[V_i][V_j]a^{k-1}}=\overline{[V_i][V_{\tau^{k-1}(j)}]}=\delta_{i,\tau^{1-k}(j^*)}\overline{[V_1]}+\sum_{t=2}^{m}\gamma_t\overline{[V_t]},$$ for some non-negative integers $\gamma_t$. Hence the coefficient $p_{(i,k)(j,k)}^{(1,1)}$ of $b_{(1,1)}$ in $b_{(i,k)}b_{(j,k)}$ is \begin{align*}p_{(i,k)(j,k)}^{(1,1)}&=\varepsilon(\overline{M[i,k]})\varepsilon(\overline{M[j,k]})\delta_{i,\tau^{1-k}(j^*)}\\
&=\varepsilon(\overline{M[i,k]})\varepsilon(\overline{M[j,k]})\delta_{(i,k),(j,k)^*}\\
&=\varepsilon(\overline{M[i,k]})\varepsilon(\overline{M[i,k]^*})\delta_{(i,k),(j,k)^*}\\
&=\varepsilon(b_{(i,k)})\delta_{(i,k),(j,k)^*}.
\end{align*}
Therefore, the condition (G3) is satisfied.
\epf

The notion of a bi-Frobenius algebra was introduced by Doi and Takeuchi in \cite{DT}. It generalizes the notion of a finite dimensional Hopf algebra. A bi-Frobenius algebra $A$ is both a Frobenius algebra and Frobenius coalgebra together with two elements $\phi\in A^*$ and $t\in A$ satisfying some special conditions. We refer to \cite{DT, Doi1, Doi3} for the definition. A group-like algebra can be viewed as a bi-Frobenius algebra in a natural way, see \cite[Example 3.2]{Doi2}. Following the similar approach, we define on $(R_{st}(H),\varepsilon,\mathbf{B},*)$ a bi-Frobenius algebra structure as follows:

$(R_{st}(H),\phi)$ is a Frobenius algebra with the Frobenius homomorphism $\phi$ given by $\phi(b_{(i,j)})=\delta_{(1,1),(i,j)}$ for $(i,j)\in I$. The bases $\{\frac{b_{(i,j)^*}}{\varepsilon(b_{(i,j)})}\mid(i,j)\in I\}$ and $\{b_{(i,j)}\mid(i,j)\in I\}$ form a dual pair of $(R_{st}(H),\phi)$. On the other hand, $R_{st}(H)$ is a coalgebra with counit $\varepsilon$, and comultiplication $\bigtriangleup$ defined by $\bigtriangleup(b_{(i,j)})=\frac{1}{\varepsilon(b_{(i,j)})}b_{(i,j)}\otimes b_{(i,j)}$ for $(i,j)\in I$. Now let $t=\sum_{(i,j)\in I}b_{(i,j)}$. Then $(R_{st}(H),t)$ is a Frobenius coalgebra. Define a map $$S:R_{st}(H)\rightarrow R_{st}(H),\ S(b_{(i,j)})=b_{(i,j)^*}\ \textrm{for}\ (i,j)\in I.$$ It is easy to see that the map $S$ is an anti-algebra and anti-coalgebra map, so called an antipode of $R_{st}(H)$. The quadruple $(R_{st}(H),\phi,t,S)$ forms a bi-Frobenius algebra since it satisfies the defining conditions of a bi-Fronbenius algebra. That is, $\varepsilon$ is an algebra morphism from $R_{st}(H)$ to $\mathbb{C}$; $b_{(1,1)}=1$ is the group-like element of $R_{st}(H)$; $(R_{st}(H),\phi)$ is a Frobenius algebra; $(R_{st}(H),t)$ is a Frobenius coalgebra; and $S$ is an antipode. As a consequence, various properties of group-like
algebras and bi-Frobenius algebras (see \cite{Doi1, Doi2, Doi3, DT}) hold for $R_{st}(H)$.

We have given $R_{st}(H)$ a bi-Frobenius structure, where the coalgebra structure of $R_{st}(H)$ is defined directly on indecomposable modules of $H$. Since the stable Green algebra $R_{st}(H)$ is isomorphic to $R(\kk G)[z]/(F_n(a,z))$, we would like to see the corresponding coalgebra structure on $R(\kk G)[z]/(F_n(a,z))$. From now on,  we will identify $R_{st}(H)$ with $R(\kk G)[z]/(F_n(a,z))$.

Recall that the Dickson polynomial $F_s(y,z)$ of the second type is defined recursively by (\ref{equ7}).
The general form of the polynomial $F_s(y,z)$ is given in \cite[Lemma3.3]{LZ}:
\be\label{equ7.1}F_s(y,z)=\sum_{i=0}^{\lfloor\frac{s-1}{2}\rfloor}(-1)^i\begin{bmatrix}
                                                            s-1-i \\
                                                            i \\
                                                          \end{bmatrix}
y^iz^{s-1-2i},\ee
where $s\geq 2$ and $\lfloor\frac{s-1}{2}\rfloor$ stands for the biggest integer not bigger than $\frac{s-1}{2}$.
It follows from Equation (\ref{equ7.1}) that $\{\overline{[V_i]z^j}\mid 1\leq i\leq m, 0\leq j\leq n-2\}$ forms a basis of $R(\kk G)[z]/(F_n(a,z))$.
We have the following inverse version of the Dickson polynomials.

\begin{lem}\label{lem7}
The following equation holds for all $s\geq 1$:
$$z^s=\sum_{i=0}^{\lfloor\frac{s}{2}\rfloor}\begin{pmatrix}
                                                           s \\
                                                           i \\
                                                         \end{pmatrix}
\frac{s+1-2i}{s+1-i}y^iF_{s+1-2i}(y,z).$$
\end{lem}
\pf By induction on $s$ one can verify that the equation holds in both cases where $s$ is odd  or even.\epf

 From the identification of $R_{st}(H)$ with $R(\kk G)[z]/(F_n(a,z))$, one can see that the basis elements $\overline{M[i,j]}$ in $R_{st}(H)$ correspond to $\overline{[V_i]F_j(a,z)}$ in $R(\kk G)[z]/(F_n(a,z))$, see Theorem \ref{th4.1} and Proposition \ref{prop6.1}. On the other hand, by Lemma \ref{lem7}, one has the following expressions for  the basis elements $\overline{[V_i]z^j}$ in $R(\kk G)[z]/(F_n(a,z))$.
\begin{align*}
\overline{[V_i]z^j}&=\sum_{k=0}^{\lfloor\frac{j}{2}\rfloor}\begin{pmatrix}
                                                           j \\
                                                           k \\
                                                         \end{pmatrix}
\frac{j+1-2k}{j+1-k}\overline{[V_i]a^kF_{j+1-2k}(a,z)}\\
&=\sum_{k=0}^{\lfloor\frac{j}{2}\rfloor}\begin{pmatrix}
                                                           j \\
                                                           k \\
                                                         \end{pmatrix}
\frac{j+1-2k}{j+1-k}\overline{[V_{\tau^k(i)}]F_{j+1-2k}(a,z)})\\
&=\sum_{k=0}^{\lfloor\frac{j}{2}\rfloor}\begin{pmatrix}
                                                           j \\
                                                           k \\
                                                         \end{pmatrix}
\frac{j+1-2k}{j+1-k}\overline{M[\tau^k(i),j+1-2k]}.
\end{align*}
Now we can summarize  the bi-Frobenius algebra structure on $R(\kk G)[z]/(F_n(a,z))$ as follows:

\begin{prop}
The quadruple $(R(\kk G)[z]/(F_n(a,z)),\phi,t,S)$ forms a bi-
Frobenius algebra. The Frobenius homomorphism $\phi$, the comultiplication $\Delta$, the element $t$ and the antipode $S$ are respectively given by

\begin{itemize}
\item $\phi(\overline{[V_i]z^j})=\begin{cases}
\begin{pmatrix}
                                                           j \\
                                                           \frac{j}{2} \\
                                                         \end{pmatrix}
\frac{2}{j+2}, & j\ \textrm{is}\ \textrm{even}\ \textrm{and}\ [V_i]=a^{-\frac{j}{2}},\\
0, & \textrm{otherwise},
\end{cases}$ 
\item 
$\bigtriangleup(\overline{[V_i]z^j})=\sum_{k=0}^{\lfloor\frac{j}{2}\rfloor}\begin{pmatrix}
                                                           j \\
                                                           k \\
                                                         \end{pmatrix}
\frac{j+1-2k}{j+1-k}\frac{1}{\dim(V_i)F_{j+1-2k}(1,2\cos\frac{\pi}{n})}$\\
$\cdot\overline{[V_{\tau^k(i)}]F_{j+1-2k}(a,z)}\otimes\overline{[V_{\tau^k(i)}]F_{j+1-2k}(a,z)},$
\item
$t=\sum_{(i,j)\in I}b_{(i,j)}=\sum_{(i,j)\in I}\dim(V_i)F_j(1,2\cos\frac{\pi}{n})\overline{[V_i]F_j(a,z)},$
\item
$S(\overline{[V_i]z^j})=\sum_{k=0}^{\lfloor\frac{j}{2}\rfloor}\begin{pmatrix}
                                                           j \\
                                                           k \\
                                                         \end{pmatrix}
\frac{j+1-2k}{j+1-k}\overline{[V_{i^*}]a^{k-j}F_{j+1-2k}(a,z)}.$
\end{itemize}
\end{prop}

In general, the complexified Green algebra of a finite dimensional pointed Hopf algebra of finite representation type is not necessarily a group-like algebra. But we don't know whether or not the complexified stable Green algebra is a bi-Frobenius algebra.

\end{document}